\documentclass[a4paper]{amsart}

\usepackage{setspace}
\usepackage{color}
\usepackage{hyperref}
\usepackage{tikz}
\usetikzlibrary{external}
\usetikzlibrary{decorations.pathreplacing}
\usetikzlibrary{backgrounds}
\usepackage{tkz-graph}
\usepackage{dirtytalk}
\usepackage[shortlabels]{enumitem}
\usepackage{bm}
\usepackage{mathtools}

\newcommand{\G}{\mathcal{G}}

\newcommand{\V}{\mathcal{V}}
\newcommand{\I}{\mathcal{I}}

\newcommand{\graphs}{\textbf{G}}
\renewcommand{\L}{\mathcal{L}}

\newcommand{\Pt}{\mathrm{P}}
\newcommand{\blank}{\ \bullet \ }

\newcommand{\bN}{\mathbb{N}}
\newcommand{\bZ}{\mathbb{Z}}
\newcommand{\bR}{\mathbb{R}}

\newcommand{\ip}{\textrm{I}}
\newcommand{\mip}{\textrm{M}}
\newcommand{\cmip}{\textrm{C}}
\newcommand{\indep}{\textrm{Indep}}
\newcommand{\stab}{\mathcal{S}}
\newcommand{\map}[1]{\operatorname{MAP}_{#1}}
\newcommand{\lab}[1]{\normalfont \textbf{#1}}

\newcommand{\card}[1]{|#1|}
\newcommand{\homo}[1]{\mathop{\textrm{H}}(#1)}
\newcommand{\red}[1]{\tilde{#1}}
\newcommand{\defeq}{\vcentcolon=}

\DeclareMathOperator{\Ima}{\normalfont \textrm{Im}}
\DeclareMathOperator{\gl}{\normalfont \textrm{Gl}}
\DeclareMathOperator{\supp}{\normalfont \textrm{Supp}}
\DeclareMathOperator{\sym}{\normalfont \textrm{sym}}
\DeclareMathOperator{\Tr}{\normalfont \textrm{Tr}}

\newenvironment{mycase}[1]
{\par\smallskip\noindent \emph{Case #1.}
  \noindent\ignorespaces }
{\par\smallskip }

\theoremstyle{plain}
\newtheorem{THM}{Theorem}[section]
\newtheorem{PROP}[THM]{Proposition}
\newtheorem{LEMMA}[THM]{Lemma}

\newtheorem{COR}[THM]{Corollary}
\theoremstyle{definition}
\newtheorem{REM}[THM]{Remark}
\newtheorem{EX}[THM]{Example}
\newtheorem{DEF}[THM]{Definition}
\newtheorem*{DEF*}{Definition}

\theoremstyle{plain}
\newtheorem*{futuretheoreminner}{Theorem~\thefuturetheoreminner}
\newcommand{\thefuturetheoreminner}{} 

\ExplSyntaxOn

\prop_new:N \g_alevel_future_prop

\NewDocumentEnvironment{futuretheorem}{ m o +b }
{\renewcommand{\thefuturetheoreminner}{\ref{#1}}
  \IfNoValueTF{#2}
  {\futuretheoreminner}
  {\futuretheoreminner[#2]}
  #3
}
{
  \endfuturetheoreminner
  \prop_gput:Nnn \g_alevel_future_prop { #1 } { #3 }
  \IfValueT{#2}{ \prop_gput:Nnn \g_alevel_future_prop { #1-attr } { #2 } }
}

\NewDocumentCommand{\pasttheorem}{m}
{
  \prop_if_in:NnTF \g_alevel_future_prop { #1-attr }
  {
    \begin{THM}[\prop_item:Nn \g_alevel_future_prop { #1-attr }]
    }
    {
      \begin{THM}
      }
      \label{#1}
      \prop_item:Nn \g_alevel_future_prop { #1 }
    \end{THM}
  }

  \ExplSyntaxOff

  \title{Lorentzian polynomials and the independence sequences of graphs}
  \author{Amire Bendjeddou, Leonard Hardiman}
  \address{CSFT/Institute of Mathematics, EPFL, Lausanne, Switzerland}
  \email{amire.bendjeddou@epfl.ch}
  \email{leonard.hardiman@epfl.ch}
  \date{\today}

  \usepackage{lipsum}


\begin{document}

  \begin{abstract}
    We study the multivariate independence polynomials of graphs and the log-concavity of the coefficients of their univariate restrictions. Let $R_{W_4}$ be the operator defined on simple and undirected graphs which replaces each edge with a caterpillar of size $4$. We prove that all graphs in the image of $R_{W_4}$ are what we call pre-Lorentzian, that is, their multivariate independence polynomial becomes Lorentzian after appropriate manipulations. In particular, as pre-Lorentzian graphs have log-concave (and therefore unimodal) independence sequences, our result makes progress on a conjecture of Alavi, Malde, Schwenk~and Erd\H{o}s which asks if the independence sequence of trees or forests is unimodal.
  \end{abstract}

  \maketitle

  \section{Introduction}
  A subset of vertices in a finite, simple, undirected graph is called \emph{independent} if it contains no two vertices that are adjacent. For such a graph $G$, one can consider the so-called \emph{independence sequence},
  \begin{align*}
    i_k = \#\{S \subset V(G) \mid S \ \text{is independent and} \ \#S = k \}.
  \end{align*}
  As is often the case, considering the corresponding generating function, the \emph{independence polynomial}, has proved valuable in studying this sequence. In this paper, we develop a novel methodology for establishing various properties of the independence polynomial and provide an application. To contextualise these results within the existing literature, we now briefly discuss related results.  

  We first consider a closely related object to the independence polynomial: the matching polynomial; indeed the matching polynomial of a graph $G$ is the independence polynomial of the line graph of $G$. It is well known that the matching polynomial is real rooted and more recently Chudnovsky and Seymour~\cite{chudnovsky2007roots} proved that this is true for the independence polynomial of any claw-free graph. Both the matching and independence polynomial have natural multivariate versions by associating each vertex to a unique variable. The multivariate version of the matching polynomial was  proven to be real stable by Heilmann
  and Lieb~\cite{MR297280}, which can be seen as a multivariate generalisation of real rootedness. In a similar spirit, Leake and Ryder~\cite{ALCO_2019__2_5_781_0} show that the Chudnovsky-Seymour result can be generalised to the multivariate independence polynomial of claw-free graphs, by introducing the notion of same-phase stability.

  Before proceeding, we take a moment to state certain key definitions that the following discussion will rely upon.
  \begin{DEF}
    \label{def:log-concave}
    A non-negative sequence $(s_0, s_1,\dots,s_n)$ is called \emph{log-concave} if
    \begin{align*}
      s_{k}^{2} \geq s_{k-1}s_{k+1}, \ \forall k \in \{1, 2, \ldots, n-1\}.
    \end{align*}
  \end{DEF}
  \begin{DEF}
    \label{def:ultra-log-concave}
    A non-negative sequence $(s_0, s_1, \dots,s_n)$ is called \emph{ultra log-concave} if
    \begin{align*}
      \frac{s_{k}^{2}}{\binom{n}{k}^{2}} \geq \frac{s_{k-1}}{\binom{n}{k-1}} \frac{s_{k+1}}{\binom{n}{k+1}}, \ \forall k \in \{1, 2, \ldots, n-1\}.
    \end{align*}
    Note that, as the sequence $s_k = \binom{n}{k}$ is log-concave, ultra log-concavity implies log-concavity.
  \end{DEF}
  It is well known that the coefficients of real rooted polynomials form ultra log-concave sequences~\cite[Theorem B, p.\ 270]{MR0460128}. Motivated by the mentioned successful approaches of considering the multivariate independence polynomial to prove real rootedness (and hence ultra-log concavity), in this paper we show that Brändén and Huh's recently developed theory of Lorentzian polynomials~\cite{MR4172622} can also be used to prove that the independence polynomial of certain families of graphs have log-concave coefficients. In particular, we introduce a property of partitioned graphs (essentially graphs with a distinguished subset of vertices, see Definition~\ref{def:partitioned-graph}) called \emph{pre-Lorentzian}, which implies log-concavity of the coefficients of their independence polynomials. Let us now introduce the following operator on graphs. 
  \begin{DEF}\label{def:replace-op}
    Let $R_{W_4}$ denote the operator on graphs which replaces every edge in a graph with the following, 
    \begin{align*}
      W_4 =
      \tikzsetnextfilename{W_3}
      \definecolor{blue}{rgb}{0.44, 0.61, 1}
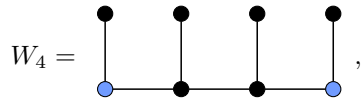
      \begin{array}{c}
        \begin{tikzpicture}[main/.style = {draw, circle, minimum size=0.2cm, inner sep=0pt, outer sep=0pt}, edge/.style = {line width = 0.2mm},scale=1,every node/.append style={transform shape}] 
          \node[main,fill=blue] (1) {};
          \node[main,fill=black, right of=1] (2) {};
          \node[main,fill=black, right of=2] (3) {};
          \node[main,fill=blue, right of=3] (4) {};
          \node[main,fill=black, above of=1] (5) {};
          \node[main,fill=black, above of=2] (6) {};
          \node[main,fill=black, above of=3] (7) {};
          \node[main,fill=black, above of=4] (8) {};
          \draw[edge] (1) -- (2);
          \draw[edge] (2) -- (3);
          \draw[edge] (3) -- (4);
          \draw[edge] (1) -- (5);
          \draw[edge] (2) -- (6);
          \draw[edge] (3) -- (7);
          \draw[edge] (4) -- (8);
        \end{tikzpicture} 
      \end{array} 
      ,
    \end{align*}
    where the blue vertices are the ones to be identified with the vertices of the original graph (we call these vertices the \emph{endpoints}).
  \end{DEF}
  To illustrate the effectiveness of our approach, we provide an application in the form of the following result.
  \begin{futuretheorem}{thm:main-theorem}
    Let $\graphs$ denote the set of simple finite undirected graphs (not necessarily connected). Then every graph in $R_{W_4}(\graphs)$ has an independence sequence that is log-concave.
  \end{futuretheorem}

  We now say a few words on the proof of Theorem~\ref{thm:main-theorem}. In truth, our principal tool is not the multivariate independence polynomial, but rather the \emph{coloured independence polynomial}, which may be thought of as a compromise between the independence polynomial and its multivariate analogue. Indeed, just as the independence polynomial may be recovered by identifying the variables in its multivariate analogue, the coloured independence polynomial is recovered by partially identifying the variables. Crucially, the coloured independence polynomial is flexible enough to satisfy the properties we require for non-trivial families of graphs (unlike the multivariate independence polynomial) and powerful enough to allow for glueing arguments (unlike the univariate independence polynomial). Indeed, these glueing arguments are at the heart of the proof of Theorem~\ref{thm:main-theorem}. In particular, we shall show that an integer indexed family of star-like graphs are pre-Lorentzian (Theorem \ref{thm:leafy-stars-lorentz}) and that any graph in the image of $R_{W_4}$ can be constructed by glueing together members of this family (Proposition \ref{prop:replace_to_glue}). 

  A pleasing feature of Theorem~\ref{thm:main-theorem} is that it implies progress on the following, well studied, conjecture. 
  Alavi, Malde, Schwenk and Erd\H{o}s \cite{alavi1987vertex} asked if the independence sequence of trees or forests is unimodal, a property implied by log-concavity. Despite substantial effort (see~\cite{MR3251932,MR3862630,MR3737109,MR2062222,MR1985180,MandrescuSpivak2016,MR2727455,MR3490439,MR3864110,MR4153624,MR2324272}, furthermore,~\cite{basit2021independent} provides an extensive survey on this topic) in proving the conjecture for different families of trees, all the constructions so far exhibit a recursive structure which restricts the families of trees that can be built. By using the coloured independence polynomial we overcome this constraint (in particular due to the possibility of the aforementioned glueing arguments). Indeed, applying our edge replacement to forests we prove the conjecture for all forests which can be built by replacing all edges of arbitrary forests with caterpillars of size 4. We note, however, that the full conjecture remains beyond the reach of the techniques developed in this paper, as demonstrated by the following result of Kadrawi and Levit. 
  \begin{THM}[\cite{KL23}]
    \label{thm:trees-not-log-concave}
    There exist infinite families of trees whose independence sequences fail to be log-concave. 
  \end{THM}

  \subsection{Structure of this paper} Section~\ref{sec:preliminaries} details the necessary background material and records certain straightforward lemmas from linear algebra. Section~\ref{sec:mip-glueing} then introduces the coloured independence polynomial, develops the theory enabling the corresponding glueing arguments, and specialises these arguments to the pre-Lorentzian property. Finally, Section~\ref{sec:application} concludes by providing an application in the form of Theorem~\ref{thm:main-theorem}.

  \subsection{Conventions} All graphs in this paper are assumed to be finite, simple and undirected. For a graph $G$, we use $V(G)$ and $E(G)$ to denote the sets of vertices and edges, respectively. 
  We will also discuss \emph{coloured} graphs in this paper, however we do not require the colouring to be proper, but instead make use of the following, simpler, concept.

  \begin{DEF}
    \label{def:coloured-graph}
    A \emph{coloured graph} $\G$ is a pair $(G,i)$, where $G$ is a graph and $i$ is a map from $\V(G)$ to some indexing set $\I$, we call $\I$ the \emph{set of colours}. 
  \end{DEF}

  Throughout the paper, we will make extensive use of the following terminology, which we record here for convenience.  

  \begin{DEF}[Free vertex]
    \label{def:free-vertex}
    For a vertex $v$ in a coloured graph, we say that $v$ (and its colour) are \emph{free} if the subset of vertices that share a colour with $v$ is $\{v\}$.
  \end{DEF}

  For $n \in \bN$, we use $[n]$ to denote the set $\{1,2,\dots,n\}$. When a collection of $n$ variables is written with a shared symbol $x$ decorated by a subscript in $[n]$ (i.e.\ $x_1,x_2, \dots, x_n$), we write, for all $\alpha \in \bN^n$, $x^{\alpha}$ to denote $\prod x^{\alpha_i}_i$. Unless otherwise specified, all polynomials are assumed to be over $\bR$. For $n \in \bN^{*}$, we use $<$ and $\leq$ to denote the partial orders on $\bN^n$ given by
  \begin{align*}
    \alpha < \beta \iff \alpha_i < \beta_i, \forall i \in [n] \quad \text{and}\quad \alpha \leq \beta \iff \alpha_i \leq \beta_i, \forall i \in [n].
  \end{align*}
  We also use the symbol $e_i$ to denote the standard basis vector $(0, \ldots, 0, 1, 0, \ldots, 0)$, where the $1$ is in the $i$th position.

  \section{Preliminaries} \label{sec:preliminaries}

  \subsection{Lorentzian polynomials}
  The notion of Lorentzian polynomials has appeared under various names throughout the literature. They were first studied  by Gurvits~\cite{MR2683227}, who referred to them as strongly log-concave polynomials, to explore multivariate generalizations of Newton’s inequalities. Years later, the theory was simultaneously advanced by two groups: Anari et al.~\cite{MR4332671, MR4681146, MR4728467} developed applications under the name completely log-concave polynomials, proving the strongest version of Mason's conjecture, while Brändén and Huh~\cite{MR4172622} coined the term \say{Lorentzian polynomial}, gave a new definition, established its equivalence with both previous formulations (in the homogeneous case) and gave an independent proof of the strongest conjecture of Mason.

  We start by recording certain properties of Lorentzian polynomials. The results in this section are well-known to experts in the field however we include their statements and proofs for the sake of completeness.

  \begin{DEF}
    \label{def:M-convex}
    A finite subset $A \subset \mathbb{Z}^n$ is called \emph{M-convex} if for any $\alpha, \beta \in A$ and any index $i$ satisfying $\alpha_i > \beta_i$, there is an index $j$ satisfying
    \begin{align*}
      \beta_j > \alpha_j \ \text{and} \   \alpha - e_i + e_j \in A.
    \end{align*}
  \end{DEF}

  \begin{DEF}
    \label{def:support}
    The \emph{support} of a multivariate polynomial $p = \sum_{\alpha \in \bN^n} c_{\alpha}x^{\alpha}$ is given by the set $\supp(p) = \{ \alpha \in \bN^n : c_{\alpha} \neq 0 \}$. 
  \end{DEF}

  \begin{DEF}
    \label{def:lorentz-polyn}
    A homogeneous polynomial $p$ of degree $d$ with non-negative coefficients is called \emph{Lorentzian} if 
    \begin{itemize}
    \item $\supp(p)$ is M-convex, 
    \item for every $(d-2)$th partial derivative of $p$, the corresponding Hessian matrix has at most one positive eigenvalue.
    \end{itemize}
  \end{DEF}

  \begin{REM}
    \label{rem:lorentz-polyn-def}
    Lorentzian polynomials were introduced by Brändén and Huh~\cite{MR4172622} using an alternative definition. In particular,~\cite[Theorem 2.25]{MR4172622} provides an equivalence between the original definition and the definition given above. 
  \end{REM}

  \begin{PROP}
    \label{prop:lorentz-polyn-stable-indent}
    Lorentzian polynomials are stable under identifying variables.
    \proof
    This is a special case of~\cite[Theorem 2.10]{MR4172622}. \endproof
  \end{PROP}

  \begin{PROP}
    \label{prop:lorentz-polyn-stable-prod}
    Lorentzian polynomials are stable under taking products.

    This was proven in ~\cite[Corollary 2.32.]{MR4172622}. For sake of completeness we outline a more detailed proof here.  
    
    \proof Following Gurvits~\cite{MR2683227}, we say that a polynomial $p(x_1,\ldots,x_n)$ in $n$ variables with non-negative coefficients is \emph{strongly log-concave} if, for any $i \in [n]$ and $a_i \in \bN$,  
    \begin{align*}
      \prod_{i=1}^n \partial_{x_i}^{a_i}p \text{ is identically zero or } \log\Big( \prod_{i=1}^n \partial_{x_i}^{a_i}p \Big) \text{ is concave on } \bR_{>0}^n.
    \end{align*}
    Homogeneous polynomials are Lorentzian if and only if  they are strongly log-concave~\cite[Theorem 2.30]{MR4172622}, so it suffices to show the result for strongly log-concave polynomials. Let $u=(u_1, \ldots, u_m)$ and $w=(w_1, \ldots, w_d)$ be two disjoint sets of variables and let $f(u)$, $g(w)$ be two homogeneous strongly log-concave polynomials. It is clear that $h(u,w)=f(u)g(w)$ is a homogeneous polynomial with non-negative coefficients. We need to show that $\prod_{i=1}^m \partial_{u_i}^{a_i} \prod_{j=1}^d  \partial_{w_j}^{b_j}h(u,w)$ is either identically zero or log-concave on $\bR_{>0}^{m+d}$, where $a_i \in \bN$ and $b_j \in \bN$ are arbitrary. We have $\prod_{i=1}^m \partial_{u_i}^{a_i} \prod_{j=1}^d  \partial_{w_j}^{b_j}h(u,w) = \prod_{i=1}^m \partial_{u_i}^{a_i} f(u) \cdot \prod_{j=1}^d\partial_{w_j}^{b_j} g(w) $. If either $\prod_{i=1}^m \partial_{u_i}^{a_i} f(u)$ or $\prod_{j=1}^d\partial_{w_j}^{b_j} g(w) $ are identically zero, we are done. If both polynomials are not identically zero, then
    \begin{align*}
      \log\Big(\prod_{i=1}^m \partial_{u_i}^{a_i} \prod_{j=1}^d  \partial_{w_j}^{b_j}h(u,w)\Big) = \log\Big(\prod_{i=1}^m \partial_{u_i}^{a_i} f(u)\Big) + \log\Big(\prod_{j=1}^d\partial_{w_j}^{b_j} g(w)\Big), 
    \end{align*}
    which is a sum of concave functions on $\bR_{>0}^{m+d}$. Since sums of concave functions are concave, $h(u,w)$ is strongly log-concave, thus Lorentzian. The general case where $u$ and $w$ share variables now follows by applying Proposition \ref{prop:lorentz-polyn-stable-indent} to $h(u,w)$. 
    \endproof
  \end{PROP}

  \begin{DEF}
    \label{def:power-truncation}
    For $\alpha, \beta \in \bN^{n}$ the \emph{power truncation operator} $\Pt_{\alpha,\beta}$ is defined as follows,
    \begin{align*}
      \Pt_{\alpha, \beta}: \quad &\bR[x_1, \ldots, x_n] \to \bR[x_1, \ldots, x_n] \\
                                 &\sum_{\gamma \in \bN^n } a(\gamma) x^{\gamma} \mapsto \sum_{\alpha \leq \gamma \leq \beta} a(\gamma)x^{\gamma}.
    \end{align*}
  \end{DEF}
  \begin{PROP}
    \label{prop:lorentz-polyn-stable-pt}
    Lorentzian polynomials are stable under $\Pt_{\alpha,\beta}$.

    \proof Let $p \in \bR[x_1,\dots,x_n]$ be a Lorentzian polynomial. By~\cite[Theorem 3.2]{MR4172622}, $\Pt_{\alpha,\beta}(p)$ being Lorentzian follows from the following two claims:
    \begin{enumerate}[(i)]
    \item there exists $l\in \bN$ such that $\Pt_{\alpha, \beta}$ is homogeneous of degree $l$, i.e. it either kills a monomial, or changes its degree by adding $l$. This holds trivially for $l=0$, as $\Pt_{\alpha,\beta}$ either kills a monomial or leaves its degree unchanged.
    \item \label{cl:hard} The polynomial 
      \begin{align*}
        \sym_{\Pt_{\alpha,\beta}}((x_i),(y_i)) \defeq \sum\limits_{0\leq \gamma \leq \kappa} \binom{\kappa}{\gamma} \Pt_{\alpha,\beta}(x^{\gamma})y^{\kappa-\gamma}
      \end{align*}
      is Lorentzian, where $\kappa \in \bN^n$  is such that $\kappa_i \geq \deg p$, and $\binom{\kappa}{\alpha} = \prod_i \binom{\kappa_i }{\alpha_i}$.
    \end{enumerate}
    To see that Claim~\ref{cl:hard} holds, observe that
    \begin{align*}
      \sum\limits_{0\leq \gamma \leq \kappa} \binom{\kappa}{\gamma} \Pt_{\alpha,\beta}(x^{\gamma})y^{\kappa-\gamma} 
      &= \sum_{\alpha \leq \gamma \leq \beta} \binom{\kappa}{\gamma}x^{\gamma}y^{\kappa - \gamma} \\
      &= \prod_{i=1}^{n} \left(\sum_{\alpha_i \leq j \leq \beta_i}\binom{k_i }{j}x_i^{j}y_i^{k_i-j}\right), 
    \end{align*}
    hence it is enough to show that 
    \begin{align*}
      \sum_{\alpha_i \leq j \leq \beta_i}\binom{k_i }{j } x_i^{j}y_i^{k_i-j}
    \end{align*} is Lorentzian for all $i \in [n]$. This follows from the fact that a bivariate homogeneous polynomial is Lorentzian if and only if the sequence of coefficients is ultra log-concave, non-negative and has no internal zeros (see~\cite[Example 2.26]{MR4172622}). 
    \endproof
  \end{PROP}

  \begin{DEF}
    \label{def:map}
    Let $R$ be a polynomial ring over $\bR$ in a finite set of variables $V$, and let $x$ be in $V$. The \emph{multi-affine part} associated to $x$ is the unique linear operator $\map{x} \colon R \to R$ which satisfies
    \begin{align*}
      \map{x}\left(\prod_{v \in V} v^{d_v}\right) = \begin{cases} \prod\limits_{v \in V} v^{d_v} & \text {if } d_x \leq 1, \\ 0 \quad \text{otherwise,} \end{cases}
    \end{align*}
    in other words, the linear operator which simply kills any monomial of degree in $x$ greater than $1$. 
  \end{DEF}

  \begin{COR}
    \label{prop:lorentz-map-stable}
    Let $R$ be a polynomial ring over $\bR$ in a finite set of variables and let $p$ be in $R$. Lorentzian polynomials are stable under $\map{x}$. 
    \proof
    After putting an order on the finite set of variables which identifies $x$ with $x_1$, this follows immediately by noting that, 
    \begin{align*}
      \Pt_{0,\beta}(p) = \map{x_1}(p),
    \end{align*}
    when $\beta = (1,N,\dots,N)$ and $N \in \bN$ is sufficiently large. The conclusion then follows from~Proposition~\ref{prop:lorentz-polyn-stable-pt}.
    \endproof
  \end{COR}

  \subsection{Linear algebra lemmas}

  In this section, we detail certain lemmas which will be used in the proof of Theorem~\ref{thm:leafy-stars-lorentz}. In particular, all of these results follow from standard arguments from linear algebra.

  \begin{LEMMA}\label{lem:lemma1}
    Let $M$ be a $(n+2)\times (n+2)$ symmetric matrix of the form 
    \begin{align*}
      M= \begin{pmatrix}
           0 & a & \cdots  & a& b &c  \\
           a& 0  & \ddots  & \vdots  & \vdots & \vdots \\
           \vdots & \ddots & \ddots & a & b & c\\
           a & \cdots & a & 0  & b &c \\ 
           b & \cdots & b  & b & d  & f\\
           c & \cdots &  c  & c& f & e\\
         \end{pmatrix}
    \end{align*}
    where $a,b,c,d,e,f$ are non-negative real numbers. Then $M$ has $n-1$ negative eigenvalues of the form $ -a$ and the corresponding eigenspace is spanned by the vectors $(e_1 - e_i)_{2 \leq i \leq n}$. Furthermore, any of the remaining $3$ eigenvalues is either equal to $-a$ or is an eigenvalue of 
    \begin{align*}
      \red{M} \defeq 
      \begin{pmatrix}
        (n-1)a& b & c \\
        nb& d & f &  \\
        nc & f & e 
      \end{pmatrix}.               
    \end{align*}
    Within the context of this lemma, we call $\red{M}$ the \emph{reduced form} of $M$.
    \proof The first assertion of the lemma simply follows by applying the matrix to the claimed eigenvectors. For the second part, note that any of the other $3$ eigenvalues is either $-a$ or has corresponding eigenvector which is orthogonal to all $(e_1 - e_i)_{2 \leq i \leq n}$, since eigenvectors of a symmetric matrix with different eigenvalues are orthogonal. This implies that these eigenvectors with, say, eigenvalue $\lambda$, have the form $v_{\lambda} = (v,v,\ldots,v,u,w)^\top$. Now it is easy to see that the eigenvalue equation $M v_{\lambda} = \lambda v_{\lambda}$ is equivalent to 
    \[
      \begin{pmatrix}
        (n-1)a& b & c \\
        nb& d & f &  \\
        nc & f & e 
      \end{pmatrix}
      \begin{pmatrix}
        v\\
        u \\
        w
      \end{pmatrix} = 
      \lambda
      \begin{pmatrix}
        v\\
        u \\
        w
      \end{pmatrix}
    \]
    i.e.\ $\lambda$ is an eigenvalue of $\red{M}$. Conversely, since $u,v,w$ are free parameters, any eigenvalue of $\red{M}$ is also an eigenvalue of $M$. \endproof  
  \end{LEMMA}

  \begin{REM}
    \label{rem:mat3-3}
    Note that $\red{M}$ has only real eigenvalues, since its eigenvalues form a subset of the eigenvalues of $M$, and $M$ is symmetric. 
  \end{REM}

  \begin{LEMMA}
    \label{lem:lemma2}
    Let $A$ be a $3 \times 3$ matrix which satisfies 
    \begin{align}
      \label{eq:assumptions}
      \Tr(A) > 0, \  A_{11} + A_{22} + A_{33} \leq 0 \ \text{and} \  \det(A) \geq 0,
    \end{align}
    where $A_{ii}$ is the minor of $A$ obtained by deleting the $i$th row and column. Then $A$ has exactly one positive eigenvalue.
    \proof The characteristic polynomial of $A$ has the form 
    \begin{align*}
      p_A(x) = x^3 - \Tr(A) x^2+ (A_{11} + A_{22} + A_{33})x - \det(A).
    \end{align*}
    By our assumptions~\eqref{eq:assumptions}, the number of sign changes in the sequence of $p_A(x)$'s coefficients is exactly $1$. By Descartes' rule of signs $p_A$ has exactly one positive root.  \endproof
  \end{LEMMA}

  \section{Multivariate independence polynomial and glueing} \label{sec:mip-glueing}

  \subsection{Independence polynomials} As has become fairly standard within the literature, we study the independence sequence via its generating function, the independence polynomial, first defined by Gutman and Harary in 1983~\cite{MR0724764}. 

  \begin{DEF}[Independence polynomial] \label{def:indep-poly}
    Let $G$ be a finite graph. The corresponding \emph{independence polynomial}, denoted $\ip(G)$, is an element of the polynomial ring $\bZ[x]$, defined by setting 
    \begin{align*}
      \ip(G\,;x) = \sum_{S \subseteq V(G)} \indep(S) \cdot x^{\card{S}},
    \end{align*}
    where $\indep(S)$ is an indicator function for whether or not $S$ is independent and $\card{S}$ denotes the cardinality of $S$.   
  \end{DEF}

  More recently, evidence has emerged that considering the analogous multivariate polynomial can be fruitful (Borcea and Brändén's simple proof that the matching polynomial is real-rooted proves a striking example~\cite{MR2534100,MR2569072}).

  \begin{DEF}[Multivariate independence polynomial] \label{def:mul-indep-poly}
    Let $G$ be a finite graph. The corresponding \emph{multivariate independence polynomial}, denoted $\mip(G)$, is an element of the polynomial ring $\mathbb{Z}[(x_v)_{v \in V(G)}]$, defined by setting 
    \begin{align*}
      \mip(G\,;(x_v)) = \sum_{S \subseteq V(G)} \indep(S) \cdot x_S,
    \end{align*}
    where $x_S = \prod_{v \in S} x_v$.  
  \end{DEF}

  We introduce what may be considered a \say{halfway house} between the independence polynomial and its multivariate counterpart: the coloured independence polynomial. Indeed, just as $\ip(G\,;x)$ may be recovered by identifying \emph{all} the variables in $\mip(G\,;(x_v))$, the coloured independence polynomial is recovered by \emph{partially} identifying the variables. 

  \begin{DEF}[Coloured independence polynomial] \label{def:lab-mul-indep-poly}
    Let $\G = (G,i)$ be a finite coloured graph. The corresponding \emph{coloured independence polynomial}, denoted $\cmip(\G)$, is an element of the polynomial ring $\mathbb{Z}[(x_l)_{l \in \Ima i}]$, defined by setting  
    \begin{align*}
      \cmip(\G\,; (x_l)) = \mip(G,(x_v))\big|_{x_v = x_{i(v)}},
    \end{align*}
    i.e.\ the polynomial obtained by identifying variables corresponding to vertices which share the same colour. We also say that a variable $x_l$ is \emph{free} if the associated colour is (see Definition~\ref{def:free-vertex}).
  \end{DEF}

  \subsection{Glueing graphs} \label{sec:glueing-arbitrary}

  We now describe a glueing procedure that takes as input two finite coloured graphs, $\G_1 = (G_1,i_1)$ and $\G_2 = (G_2,i_2)$, with colour sets $\I_1$ and $\I_2$. We may assume, up to changing the colours, that $\I_1 \cap \I_2 = \emptyset$. The procedure also takes as input two colours, $c_1, c_2 \in \Ima i_1 \sqcup \Ima i_2$. The resulting glued graph, denoted $\gl(\G_1,\G_2\,; c_1,c_2)$, is constructed as follows:
  \begin{enumerate}
  \item Let $\I_3$ be the disjoint union of $\I_1$ and $\I_2$ with $c_1$ and $c_2$ identified into a single colour which, by convention, we call $c_1$. We start by considering the disjoint union of $\G_1$ and $\G_2$, with the induced colouring on $\I_3$.
  \item We then add any missing edges to the subgraph induced by $i^{-1}(c_1) \subseteq \V(G_1) \sqcup \V(G_2)$ until it forms a clique. The result is $\gl(\G_1,\G_2\,; c_1,c_2)$.
  \end{enumerate}

  \begin{EX}
    We consider the following two coloured graphs:
    \begin{align*}
      \G_1 =
      \tikzsetnextfilename{general_glueing}
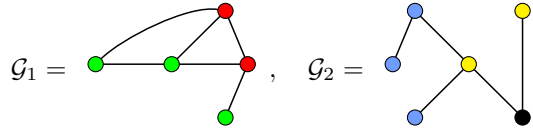
      \begin{array}{c}
        \begin{tikzpicture}[main/.style = {draw, circle, minimum size=0.2cm, inner sep=0pt, outer sep=0pt}, edge/.style = {line width = 0.2mm}]
          \node[main,fill=green] (1) {};
          \node[main,fill=red,above right of=1] (2) {};
          \node[main,fill=red,right of=1] (3) {};
          \node[main,fill=green,left of=1] (4) {};
          \node[main,fill=green,below right of=1] (5) {};
          \draw[edge] (1) -- (2);
          \draw[edge] (2) -- (3);
          \draw[edge] (1) -- (3);
          \draw[edge] (3) -- (5);
          \draw[edge] (4) -- (1);
          \draw[edge] (4) to[out=45,in=170,looseness=0.6] (2);
        \end{tikzpicture} 
      \end{array}
      , \quad
      \G_2 = \tikzsetnextfilename{general_glueing2}
      \definecolor{blue}{rgb}{0.44, 0.61, 1}
      \begin{array}{c}
        \begin{tikzpicture}[main/.style = {draw, circle, minimum size=0.2cm, inner sep=0pt, outer sep=0pt}, edge/.style = {line width = 0.2mm}] 
          \node[main,fill=yellow] (6) {};
          \node[main,fill=blue,left of=6] (7) {};
          \node[main,fill=blue,above left of=6] (8) {};
          \node[main,fill=blue,below left of=6] (9) {};
          \node[main,fill=yellow,above right of=6] (10) {};
          \node[main,fill=black,below right of=6] (11) {};
          \draw[edge] (6) -- (9);
          \draw[edge] (6) -- (8);
          \draw[edge] (11) -- (10);
          \draw[edge] (7) -- (8);
          \draw[edge] (6) -- (11);
        \end{tikzpicture} 
      \end{array}
    \end{align*}
  \end{EX}

  Then the coloured graph resulting from glueing with respect to the colours \lab{red} and \lab{blue} is given by
  \begin{align*}
    \gl(\G_1,\G_2\,; \lab{red},\lab{blue}) =
    \tikzsetnextfilename{general_glueing3}
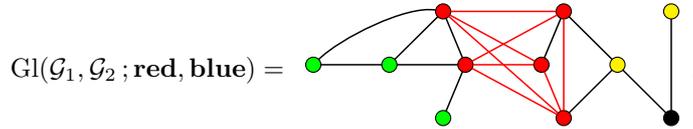
    \begin{array}{c}
      \begin{tikzpicture}[main/.style = {draw, circle, minimum size=0.2cm, inner sep=0pt, outer sep=0pt}, edge/.style = {line width = 0.2mm}] 
        \node[main,fill=green] (1) {};
        \node[main,fill=red,above right of=1] (2) {};
        \node[main,fill=red,right of=1] (3) {};
        \node[main,fill=green,left of=1] (4) {};
        \node[main,fill=green,below right of=1] (5) {};
        \draw[edge] (1) -- (2);
        \draw[edge] (2) -- (3);
        \draw[edge] (1) -- (3);
        \draw[edge] (3) -- (5);
        \draw[edge] (4) -- (1);
        \draw[edge] (4) to[out=45,in=170,looseness=0.6] (2);
        \node[main,fill=yellow,right of=3,node distance={2cm}] (6) {};
        \node[main,fill=red,left of=6] (7) {};
        \node[main,fill=red,above left of=6] (8) {};
        \node[main,fill=red,below left of=6] (9) {};
        \node[main,fill=yellow,above right of=6] (10) {};
        \node[main,fill=black,below right of=6] (11) {};
        \draw[edge] (6) -- (9);
        \draw[edge] (6) -- (8);
        \draw[edge] (11) -- (10);
        \draw[edge] (7) -- (8);
        \draw[edge] (6) -- (11);
        \draw[edge, red] (2) -- (7);
        \draw[edge, red] (2) -- (8);
        \draw[edge, red] (2) -- (9);
        \draw[edge, red] (3) -- (7);
        \draw[edge, red] (3) -- (8);
        \draw[edge, red] (3) -- (9);
        \draw[edge, red] (8) -- (9);
        \draw[edge, red] (7) -- (9);
      \end{tikzpicture}  
    \end{array}
    ,
  \end{align*}
  where the added edges have been drawn in red.

  \begin{REM}
    \label{rem:glue-on-free-vertices}
    If $i^{-1}(c_1)$ has the form $\{v_1,v_2\}$, then the only new edge after glueing is a single edge from $v_1$ to $v_2$. Note that such a glueing is possible \emph{if and only if $v_1$ and $v_2$ are free vertices\footnote{See Definition~\ref{def:free-vertex}.}} (when considered as vertices of the original unglued graphs).
  \end{REM}

  \begin{DEF}
    \label{def:glueable-stability}
    Let $\stab$ denote a property of multivariate polynomials. We call $\stab$ \emph{glueable} if the following three stability conditions hold:
    \begin{itemize}
    \item $\stab$ is stable under identifying variables,
    \item $\stab$ is stable under taking products,
    \item $\stab$ is stable under $\map{x}$ for arbitrary $x$.
    \end{itemize}
    This terminology is justified by the following proposition.
  \end{DEF}

  \begin{PROP}
    \label{prop:glueing-works}
    Let $\G_1$ and $\G_2$ be two coloured graphs such that $\cmip(\G_1)$ and $\cmip(\G_2)$ both satisfy a glueable property $\stab$. We also take $c_1$ and $c_2$ in the set of colours of $\G_1$ and $\G_2$ respectively, and consider $\G_3 = \gl(\G_1,\G_2\,; c_1,c_2)$. Then $\cmip(\G_3)$ also satisfies $\stab$.
    \proof The proposition follows from the claim that 
    \begin{align}
      \cmip(\G_3) = \map{x_{c_1}}\big(\cmip(\G_1)\cdot \cmip(\G_2)|_{c_1=c_2}\big). \label{eq:glueing-works-claim}
    \end{align}
    However, the method of constructing $\gl(\G_1,\G_2\,; c_1,c_2)$ described above makes it clear that~\eqref{eq:glueing-works-claim} holds. Indeed, $\cmip(\G_1)\cdot \cmip(\G_2) = \cmip(\G_1 \sqcup \G_2)$, and adding exactly the missing edges such that $i^{-1}(c_1)$ forms a clique has the precise effect of removing any terms in $\cmip(\G_1)\cdot \cmip(\G_2)$ which contain a power of $x_{c_1}$ greater than 1. \endproof
  \end{PROP}

  \subsection{Glueing across free vertices} \label{sec:free-glueing}

  Our aim now is to focus on glueings across free vertices. In particular, we might as well identify any colours that are not free, leading to the following definition. 

  \begin{DEF}
    \label{def:partitioned-graph}
    A \emph{partitioned graph} is a coloured graph with a distinguished colour called the \emph{bound colour} such that every colour other than the bound colour is free. As before, we extend this terminology to vertices and variables in the obvious way. 
  \end{DEF}

  \begin{REM}
    \label{rem:glueing-partitioned-trees}
    Strictly speaking, glueing two partitioned graphs across free vertices does not result in a partitioned graph, as the result will contain multiple non-free vertices. However, identifying all the non-free colours allows us to recover a partitioned graph. For the remainder of this paper, this is what is meant by the glueing of two partitioned graphs.
  \end{REM}

  An immediate corollary of Proposition~\ref{prop:glueing-works} is that two partitioned graphs whose coloured independence polynomials satisfy a glueable property $\stab$, can be glued across a free vertex into a new partitioned graph that also satisfies $\stab$. However, as we can guarantee that we will never glue across bound variables, we can in fact modify $\stab$ into a weaker property which will still be preserved by glueing as follows.

  Let $\G_1$ and $\G_2$ be partitioned graphs with free colours $c_1$ and $c_2$, respectively. By abuse of notation we use $x$ to denote the bound variable for both $\G_1$ and $\G_2$. Let us suppose that there exist $k_1,k_2 \in \bN$ such that $x^{k_1}\cdot \cmip(\G_1)$ and $x^{k_2}\cdot \cmip(\G_1)$ satisfy $\stab$, where $\stab$ is a glueable property. We now consider $\G_3$, the glueing of $\G_1$ and $\G_2$ across $c_1$ and $c_2$. As from the proof of Proposition~\ref{prop:glueing-works}, we have
  \begin{align}
    \cmip(\G_3) = \map{x_{c_1}}\big(\cmip(\G_1)\cdot \cmip(\G_2)|_{c_1=c_2}\big). 
  \end{align}
  Note that $x \notin \{x_{c_1},x_{c_2}\}$ as $c_1$ and $c_2$ are free, and therefore multiplication by $x$ commutes with $\map{c_1}(\blank|_{c_1=c_2})$ giving us
  \begin{align}
    \map{x_{c_1}}\big(x^{k_1}\cmip(\G_1)\cdot x^{k_2}\cmip(\G_2)|_{c_1=c_2}\big) = x^{k_1+k_2}\cdot\cmip(\G_3),
  \end{align}
  i.e. $x^{k_3}\cdot \cmip(\G_3)$ also satisfies $\stab$, for $k_3 = k_1 + k_2$. This discussion can be summarised into the following remark (which will play an important role in the remainder of this paper).

  \begin{REM}
    \label{rem:tree-glueing-works}
    For a glueable property $\stab$, the property of there existing $k \in \bN$ such that $x^k\cdot\cmip(\G)$ satisfies $\stab$, is stable with respect to glueing partitioned graphs across free vertices. Note that this new property is a property of the partitioned graph, not the coloured independence polynomial. 
  \end{REM}

  \subsection{Pre-Lorentzian graphs} Although there are many glueable properties worthy of study in light of Proposition~\ref{prop:glueing-works}, we choose here to focus in on one which is closely related to Lorentzian polynomials. In this section, we introduce this property.

  \begin{DEF}
    \label{def:homogenisation}
    Let $p$ be a polynomial in  $x_1, \dots, x_n$. We write $\homo{p}$ to denote the homogeneous polynomial in $x_1, \dots, x_n, y$ of minimal degree which satisfies $\homo{p}\big|_{y=1}=p$. We call $\homo{p}$ the \emph{homogenisation} of $p$ and $y$ the \emph{homogenising variable}.
  \end{DEF}

  \begin{DEF}
    \label{def:pre-lorentzian-tree}
    Let $\G$ be a partitioned graph. We call $\G$ \emph{pre-Lorentzian} if there exists $k \in \mathbb{N}$ such that $(xy)^k \cdot \homo{\cmip(\G)}$ is Lorentzian, where $x$ is the variable associated to the bound colour and $y$ is the homogenising variable.
    
    Recalling the two conditions that constitute the definition of a Lorentzian polynomial (see Definition~\ref{def:lorentz-polyn}), we note that the $(xy)^k$ factor has no impact on the M-convexity condition as a set $M \subset \bZ^n$ is M-convex if and only if $M+\alpha$ is, for all $\alpha \in \bZ^n$. It would therefore be equivalent to declare $\G$ pre-Lorentzian if  
    \begin{itemize}
    \item $\supp\big(\homo{\cmip(\G)}\big)$ is M-convex,
    \item for every $(n+2k-2)$th partial derivative of $(xy)^k \cdot \homo{\cmip(\G)}$, where $n$ is the degree of $\cmip(\G)$, the corresponding Hessian matrix has at most one positive eigenvalue.
    \end{itemize}
  \end{DEF}

  \begin{REM}
    \label{rem:pre-lorentz}
    As monomials with positive coefficients are Lorentzian, this definition is equivalent to requiring that there exist $k_1,k_2 \in \bN$ such that $x^{k_1}y^{k_2} \cdot \cmip(\G)$ is Lorentzian. 
  \end{REM}

  Due to the homogenising variable $y$ being uninvolved in the glueing procedure, we can strengthen Remark~\ref{rem:tree-glueing-works} into the following proposition.

  \begin{PROP}
    \label{prop:pre-lorentzian-glue}
    Pre-Lorentzian graphs are stable under glueing across free vertices.
    \proof
    Let $\stab$ be the property of having Lorentzian homogenisation. By Propositions~\ref{prop:lorentz-polyn-stable-indent},~\ref{prop:lorentz-polyn-stable-prod} and~\ref{prop:lorentz-polyn-stable-pt}, $\stab$ is a glueable property. The statement then follows from the exact same line of reasoning as was used to arrive at Remark~\ref{rem:tree-glueing-works} in Section~\ref{sec:free-glueing} (bolstered with the fact that multiplication by $y$ also commutes with $\map{c_1}(\blank|_{c_1=c_2})$). \endproof 
  \end{PROP}

  \begin{LEMMA}
    \label{lem:pre-lorentzian-log-concave}
    If a partitioned graph is pre-Lorentzian, then its independence sequence is log-concave.
    \proof
    Let $\G=(G,i)$ be a pre-Lorentzian partitioned graph. By identifying all free variables in $\homo{\cmip(\G)}$ with the bound variable $x$, we recover the homogenisation of the independence polynomial of $G$, i.e.,
    \begin{align*}
      \homo{\cmip(\G)}|_{x_v = x} = \homo{\ip(G;x)}.
    \end{align*}
    It follows from Proposition \ref{prop:lorentz-polyn-stable-indent}, that there exists a $k \in \bN$ such that $(xy)^k \cdot \homo{\ip(G;x)}$ is Lorentzian. Furthermore, let $d$ be the degree of $\homo{\ip(G;x)}$ and $(a_m)_{0 \leq m \leq d}$ the independence sequence of $G$, then 
    \begin{align*}
      (xy)^k \cdot \homo{\ip(G;x)} & = (xy)^k \sum_{S \subseteq V(G)} \indep(S)x^{|S|}y^{d-|S|} \\ 
                                   & = (xy)^k \sum_{m=0}^d a_m x^{m}y^{d-m} \\ 
                                   & = \sum_{m=0}^d a_m x^{m+k}y^{d+k-m}
    \end{align*}
    \noindent is a bivariate Lorentzian polynomial of degree $d+2k$ with coefficients $b_m = 0$ for $0 \leq m < k$, $b_m = a_{m-k}$ for $k \leq m < d+k$ and $b_m =0$ for $ d+k \leq m \leq d+2k$. According to~\cite[Example 2.26]{MR4172622} the coefficients $b_m$ form an ultra log-concave sequence, and since ultra log-concavity implies log-concavity, we conclude that $(a_m)_{0 \leq m \leq d}$ form a log-concave sequence. 
    \endproof 
  \end{LEMMA}

  \section{Application} \label{sec:application}

  \subsection{Characterising $R_{W_4}$ in terms of glueing} Our announced application of the material developed in Section~\ref{sec:mip-glueing} (Theorem~\ref{thm:main-theorem}) refers to an edge replacement operator whereas the material itself discusses (constrained) glueings. We now provide the bridge between these concepts.

  \begin{DEF}
    \label{def:leafy-stars}
    For $n \in \bN^{*}$, we consider the following partitioned graph on $3n+1$ vertices: 
    \begin{align*}
      \L_n =
      \tikzsetnextfilename{leafy-star}
      \definecolor{blue}{rgb}{0.44, 0.61, 1}
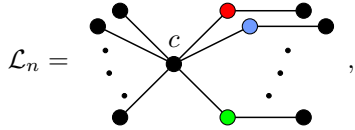
      \begin{array}{c}
        \begin{tikzpicture}[main/.style = {draw, circle, minimum size=0.2cm, inner sep=0pt, outer sep=0pt}, edge/.style = {line width = 0.2mm},scale=1,every node/.append style={transform shape}] 
          \node[main,fill=black, label={$c$}] (1) {};
          \node[main,fill=black, above left of=1] (2) {};
          \node[main,fill=black, left of=1,yshift=0.5cm] (3) {};
          \node[main,fill=black, below left of=1] (4) {};
          \node[main,fill=red, above right of=1, label={}] (5) {};
          \node[main,fill=blue, right of=1, label={},yshift=0.5cm] (6) {};
          \node[main,fill=green, below right of=1, label={}] (7) {};
          \node[main,fill=black, right of=5] (8) {};
          \node[main,fill=black, right of=6] (9) {};
          \node[main,fill=black, right of=7] (10) {};
          \node[main, minimum size=0.06cm, fill=black,yshift=0.43cm, xshift=-1cm] (12) at (2.5,-.25)  {};
          \node[main, minimum size=0.06cm, fill=black,yshift=0.43cm, xshift=-1cm] (12) at (2.4,-.55)  {};
          \node[main, minimum size=0.06cm, fill=black,yshift=0.43cm, xshift=-1cm] (12) at (2.25,-.85)  {};
          \node[main, minimum size=0.06cm, fill=black,yshift=0.43cm, xshift=-3.4cm] (12) at (2.5,-.25)  {};
          \node[main, minimum size=0.06cm, fill=black,yshift=0.43cm, xshift=-3.25cm] (12) at (2.4,-.55)  {};
          \node[main, minimum size=0.06cm, fill=black,yshift=0.43cm, xshift=-2.9cm] (12) at (2.25,-.85)  {};    
          \draw[edge] (1) -- (2);
          \draw[edge] (1) -- (3);
          \draw[edge] (1) -- (4);
          \draw[edge] (1) -- (5);
          \draw[edge] (1) -- (6);
          \draw[edge] (1) -- (7);
          \draw[edge] (5) -- (8);
          \draw[edge] (6) -- (9);
          \draw[edge] (7) -- (10);
        \end{tikzpicture} 
      \end{array} 
      ,
    \end{align*}
    where the bound vertices have been drawn in black and the vertex labelled $c$ is called the centre.
  \end{DEF}

  \begin{PROP}
    \label{prop:replace_to_glue}
    Every connected graph in $R_{W_4}(\graphs)$ may be constructed by glueing copies of $\L_n$ across free vertices.
    \proof Let $G$ be a connected graph. For each vertex $v_i$ in $V(G)$, let $\L^{(i)}_{\deg(v_i)}$ denote a copy of $\L_{\deg(v_i)}$ and let $c^{(i)}$ denote its centre. Glue $\L^{(i)}_{\deg(v_i)}$ and $\L^{(j)}_{\deg(v_j)}$ at any of the free vertices if and only if $\{v_i, v_j\} \in E(G)$. We claim that the graph obtained in this way is exactly $R_{W_4}(G)$.
    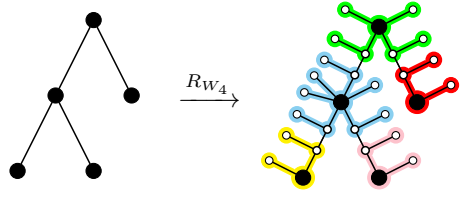
\begin{figure}[ht]
      \begin{align*}
        \begin{array}{c}
          \tikzsetnextfilename{edge-replaceA1}
          \begin{tikzpicture}[main/.style = {draw, circle, fill=black, minimum size=0.2cm, inner sep=0pt, outer sep=0pt}, edge/.style = {draw, line width = 0.2mm}, scale=1, every node/.append style={transform shape}]  
            \node[main] (1) at (0,0) {};
            \node[main] (2) at (-0.5,-1) {}; \node[main] (3) at (0.5,-1) {};  
            \node[main] (4) at (-1,-2) {}; \node[main] (5) at (0,-2) {}; 
            \draw[edge] (1) -- (2);
            \draw[edge] (1) -- (3);
            \draw[edge] (2) -- (4);
            \draw[edge] (2) -- (5);
          \end{tikzpicture}
        \end{array}
        \ \
        \xrightarrow{R_{W_4}}
        \begin{array}{c}
          \tikzsetnextfilename{edge-replaceA2}
          \definecolor{babyblue}{rgb}{0.54, 0.81, 0.94}
          \definecolor{bubblegum}{rgb}{0.99, 0.76, 0.8}
          \begin{tikzpicture}[main/.style = {draw, circle, fill=black, minimum size=0.2cm, inner sep=0pt, outer sep=0pt}, invisiedge/.style = {line width=0}, inter/.style = {draw, circle, fill=white, minimum size=0.1cm, inner sep=0pt, outer sep=0pt}, leaf/.style = {draw, circle, fill=white, minimum size=0.1cm, inner sep=0pt, outer sep=0pt}, edge/.style = {draw, line width = 0.2mm}, scale=1, every node/.append style={transform shape}]
            \node[main] (1) at (0,0) {};           
            \node[main] (2) at (-0.5,-1) {}; \node[main] (3) at (0.5,-1) {};  
            \node[main] (4) at (-1,-2) {}; \node[main] (5) at (0,-2) {}; 
            \draw[edge] (1) -- (2) node[inter, pos=0.33] (i1) {} node[inter, pos=0.67] (i2) {};
            \draw[edge] (1) -- (3) node[inter, pos=0.33] (i3) {} node[inter, pos=0.67] (i4) {};
            \draw[edge] (2) -- (4) node[inter, pos=0.33] (i5) {} node[inter, pos=0.67] (i6) {};
            \draw[edge] (2) -- (5) node[inter, pos=0.33] (i7) {} node[inter, pos=0.67] (i8) {};            
            \foreach \parent/\xshift/\yshift/\name in {i1/-0.4/0.2/l1, i2/-0.4/0.2/l2, i3/0.4/0.2/l3, i4/0.4/0.2/l4, i5/-0.4/0.2/l5, i6/-0.4/0.2/l6, i7/0.4/0.2/l7, i8/0.4/0.2/l8} {
              \node[inter] (\name) at ([xshift=\xshift cm,yshift=\yshift cm]\parent) {};
            }            
            \foreach \source/\dest in {i1/l1, i2/l2, i3/l3, i4/l4, i5/l5, i6/l6, i7/l7, i8/l8} {
              \draw[edge] (\source) -- (\dest);
            }            
            \foreach \parent/\angle/\dist/\name in {1/153/0.5/m1, 1/27/0.5/m2, 2/135/0.5/m3, 2/165/0.5/m4, 2/27/0.5/m5, 3/27/0.5/m7, 4/153/0.5/m9, 5/27/0.5/m10} {
              \node[leaf] (\name) at ([shift={(\angle:\dist)}]\parent) {};
              \draw[edge] (\parent) -- (\name);
            }
            \begin{scope}[on background layer]
              \tikzset{main/.style={draw, circle, fill=black, minimum size=0.3cm, inner sep=0pt, outer sep=0pt}}
              \tikzset{inter/.style={draw, circle, fill=blue, minimum size=0.2cm, inner sep=0pt, outer sep=0pt, line width = 0.2mm}}
              \tikzset{leaf/.style={draw, circle, fill=blue, minimum size=0.2cm, inner sep=0pt, outer sep=-1pt}}
              \tikzset{edge/.style={draw, line width = 1mm}}
              \node[main, green] (1) at (0,0) {};              
              \node[main, babyblue] (2) at (-0.5,-1) {}; \node[main, red] (3) at (0.5,-1) {};                
              \node[main, yellow] (4) at (-1,-2) {}; \node[main, bubblegum] (5) at (0,-2) {}; 
              \draw[invisiedge] (1) -- (2) node[inter, pos=0.31,green] (i1) {} node[inter, pos=0.69,babyblue] (i2) {};
              \draw[invisiedge] (1) -- (3) node[inter, pos=0.31,green] (i3) {} node[inter, pos=0.69,red] (i4) {};
              \draw[invisiedge] (2) -- (4) node[inter, pos=0.31,babyblue] (i5) {} node[inter, pos=0.69,yellow] (i6) {};
              \draw[invisiedge] (2) -- (5) node[inter, pos=0.31,babyblue] (i7) {} node[inter, pos=0.69,bubblegum] (i8) {};              
              \foreach \parent/\xshift/\yshift/\name/\parentcolor in {i1/-0.4/0.2/l1/green, i2/-0.4/0.2/l2/babyblue, i3/0.4/0.2/l3/green, i4/0.4/0.2/l4/red, i5/-0.4/0.2/l5/babyblue, i6/-0.4/0.2/l6/yellow, i7/0.4/0.2/l7/babyblue, i8/0.4/0.2/l8/bubblegum} {
                \node[inter, \parentcolor] (\name) at ([xshift=\xshift cm,yshift=\yshift cm]\parent) {};
              }              
              \draw[edge, green] (i1) -- (l1);
              \draw[edge, babyblue] (i2) -- (l2);
              \draw[edge, green] (i3) -- (l3);
              \draw[edge, red] (i4) -- (l4);
              \draw[edge, babyblue] (i5) -- (l5);
              \draw[edge, yellow] (i6) -- (l6);
              \draw[edge, babyblue] (i7) -- (l7);
              \draw[edge, bubblegum] (i8) -- (l8);              
              \foreach \parent/\angle/\dist/\name/\parentcolor in {1/153/0.5/m1/green, 1/27/0.5/m2/green, 2/135/0.5/m3/babyblue, 2/165/0.5/m4/babyblue, 2/27/0.5/m5/babyblue, 3/27/0.5/m7/red, 4/153/0.5/m9/yellow, 5/27/0.5/m10/bubblegum} {
                \node[leaf, \parentcolor] (\name) at ([shift={(\angle:\dist)}]\parent) {}; \draw[edge, \parentcolor] (\parent) -- (\name);
              }
              \draw[edge, green] (1) -- (i3);
              \draw[edge, green] (1) -- (i1);
              \draw[edge, red] (3) -- (i4);
              \draw[edge, babyblue] (i5) -- (i2);
              \draw[edge, babyblue] (2) -- (i7);
              \draw[edge, yellow] (4) -- (i6);
              \draw[edge, bubblegum] (5) -- (i8);
            \end{scope}
          \end{tikzpicture}
        \end{array}
      \end{align*}
      \caption{Glueing an edge replaced tree.}
      \label{fig:edge-replacementA}
    \end{figure}
    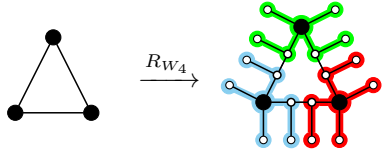
\begin{figure}[ht]
      \begin{align*}
        \begin{array}{c}
          \tikzsetnextfilename{edge-replaceB1}
          \begin{tikzpicture}[main/.style = {draw, circle, fill=black, minimum size=0.2cm, inner sep=0pt, outer sep=0pt}, edge/.style = {draw, line width = 0.2mm}, scale=1, every node/.append style={transform shape}]            
            \node[main] (1) at (0,0) {};            
            \node[main] (2) at (-0.5,-1) {}; \node[main] (3) at (0.5,-1) {};              
            \draw[edge] (1) -- (2);
            \draw[edge] (1) -- (3);
            \draw[edge] (2) -- (3);
          \end{tikzpicture}
        \end{array}
        \ \
        \xrightarrow{R_{W_4}}
        \begin{array}{c}
          \tikzsetnextfilename{edge-replaceB2}
          \definecolor{babyblue}{rgb}{0.54, 0.81, 0.94}
          \definecolor{bubblegum}{rgb}{0.99, 0.76, 0.8}
          \begin{tikzpicture}[main/.style = {draw, circle, fill=black, minimum size=0.2cm, inner sep=0pt, outer sep=0pt}, invisiedge/.style = {line width=0}, inter/.style = {draw, circle, fill=white, minimum size=0.1cm, inner sep=0pt, outer sep=0pt}, leaf/.style = {draw, circle, fill=white, minimum size=0.1cm, inner sep=0pt, outer sep=0pt}, edge/.style = {draw, line width = 0.2mm}, scale=1, every node/.append style={transform shape}]
            \node[main] (1) at (0,0) {};            
            \node[main] (2) at (-0.5,-1) {}; \node[main] (3) at (0.5,-1) {};  
            \draw[edge] (1) -- (2) node[inter, pos=0.33] (i1) {} node[inter, pos=0.67] (i2) {};
            \draw[edge] (1) -- (3) node[inter, pos=0.33] (i3) {} node[inter, pos=0.67] (i4) {};
            \draw[edge] (2) -- (3) node[inter, pos=0.33] (i5) {} node[inter, pos=0.67] (i6) {};            
            \foreach \parent/\xshift/\yshift/\name in {i1/-0.4/0.2/l1, i2/-0.4/0.2/l2, i3/0.4/0.2/l3, i4/0.4/0.2/l4, i5/0/-0.5/l5, i6/0/-0.5/l6} {
              \node[inter] (\name) at ([xshift=\xshift cm,yshift=\yshift cm]\parent) {};
            }            
            \foreach \source/\dest in {i1/l1, i2/l2, i3/l3, i4/l4, i5/l5, i6/l6} {
              \draw[edge] (\source) -- (\dest);
            }            
            \foreach \parent/\angle/\dist/\name in {1/153/0.5/m1, 1/27/0.5/m2, 2/153/0.5/m3, 2/270/0.5/m4, 3/27/0.5/m5, 3/270/0.5/m6} {
              \node[leaf] (\name) at ([shift={(\angle:\dist)}]\parent) {};
              \draw[edge] (\parent) -- (\name);
            }
            \begin{scope}[on background layer]
              \tikzset{main/.style={draw, circle, fill=black, minimum size=0.3cm, inner sep=0pt, outer sep=0pt}}
              \tikzset{inter/.style={draw, circle, fill=blue, minimum size=0.2cm, inner sep=0pt, outer sep=0pt, line width = 0.2mm}}
              \tikzset{leaf/.style={draw, circle, fill=blue, minimum size=0.2cm, inner sep=0pt, outer sep=-1pt}}
              \tikzset{edge/.style={draw, line width = 1mm}}
              \node[main, green] (1) at (0,0) {};              
              \node[main, babyblue] (2) at (-0.5,-1) {}; \node[main, red] (3) at (0.5,-1) {};  
              \draw[invisiedge] (1) -- (2) node[inter, pos=0.31,green] (i1) {} node[inter, pos=0.69,babyblue] (i2) {};
              \draw[invisiedge] (1) -- (3) node[inter, pos=0.31,green] (i3) {} node[inter, pos=0.69,red] (i4) {};
              \draw[invisiedge] (2) -- (3) node[inter, pos=0.31,babyblue] (i5) {} node[inter, pos=0.69,red] (i6) {};                            
              \foreach \parent/\xshift/\yshift/\name/\parentcolor in {i1/-0.4/0.2/l1/green, i2/-0.4/0.2/l2/babyblue, i3/0.4/0.2/l3/green, i4/0.4/0.2/l4/red, i5/0/-0.5/l5/babyblue, i6/0/-0.5/l6/red} {
                \node[inter, \parentcolor] (\name) at ([xshift=\xshift cm,yshift=\yshift cm]\parent) {};
              }              
              \draw[edge, green] (i1) -- (l1);
              \draw[edge, babyblue] (i2) -- (l2);
              \draw[edge, green] (i3) -- (l3);
              \draw[edge, red] (i4) -- (l4);
              \draw[edge, babyblue] (i5) -- (l5);
              \draw[edge, red] (i6) -- (l6);
              \foreach \parent/\angle/\dist/\name/\parentcolor in {1/153/0.5/m1/green, 1/27/0.5/m2/green, 2/155/0.5/m3/babyblue, 2/270/0.5/m4/babyblue, 3/27/0.5/m5/red, 3/270/0.5/m6/red} {
                \node[leaf, \parentcolor] (\name) at ([shift={(\angle:\dist)}]\parent) {}; \draw[edge, \parentcolor] (\parent) -- (\name);
              }
              \draw[edge, green] (1) -- (i3);
              \draw[edge, green] (1) -- (i1);
              \draw[edge, red] (3) -- (i4);
              \draw[edge, red] (3) -- (i6);
              \draw[edge, babyblue] (i5) -- (2);
              \draw[edge, babyblue] (i2) -- (2);
            \end{scope}
          \end{tikzpicture}
        \end{array}
      \end{align*}
      \caption{Glueing an edge replaced cycle.}
      \label{fig:edge-replacementB}
    \end{figure}
    For $\{v_i,v_j\} \in E(G)$, let $e_{ij}$ be the unique edge connecting $\L^{(i)}_{\deg(v_i)}$ and $\L^{(j)}_{\deg(v_j)}$. It is clear that, for every such edge, we can choose a unique $W_4$, that contains $e_{ij}$ with endpoints (see Definition~\ref{def:replace-op}) $c^{(i)}$ and $c^{(j)}$. In particular, this $W_4$ can be seen as the image of $e_{ij}$ under $R_{W_4}$. \endproof
    
  \end{PROP}

  \subsection{Proof of the main result} With Proposition~\ref{lem:pre-lorentzian-log-concave} and~\ref{prop:replace_to_glue} in mind, the main obstacle to proving Theorem~\ref{thm:main-theorem} becomes clear: establishing that $\L_n$ is pre-Lorentzian. We now proceed with the proof.

  \begin{THM}
    \label{thm:leafy-stars-lorentz}
    $\L_n$ is pre-Lorentzian for all $n \in \bN$.
    \proof
    We start by noting that the coloured independence polynomial of $\L_n$ is given by
    \begin{align}
      \begin{split}\label{eq:Ln-indep-poly}
        \cmip(\L_n;x,(x_i)) &=  x(x+1)^n + (x+1)^n \prod_{i=1}^{n}(x+x_i+1)\\
                            &=(x+1)^n \left(x + \prod_{i=1}^{n}(x+x_i+1)\right)    
      \end{split}
    \end{align}
    where $x$ is the variable associated to the bound colour and the $x_i$ are the variables associated to the free colours. To see that~\eqref{eq:Ln-indep-poly} holds, consider the fact that an independent set can either contain the central vertex $c$, giving the $x(x+1)^n$ term, or not, giving the $(x+1)^n \prod_{i=1}^{n}(x+x_i+1)$ term.

    We recall that to prove that $\L_n$ is pre-Lorentzian, we have to show that the homogenisation
    \begin{align*}
      \homo{\cmip(\L_n;x,(x_i))} = (x+y)^n \left(xy^{n-1} + \prod_{i=1}^{n}(x+x_i+y)\right)
    \end{align*}
    is Lorentzian after being multiplied by $(xy)^k$, for some $k \in \bN$. In proving this, it will be useful to suppose $n>2$, we therefore deal with the $n=1,2$ cases now.

    For $n=1$, $\homo {L(\L_1, x, (x_1)) }= (x+y)(2x+x_1+y)$ is a product of Lorentzian polynomials, thus Lorentzian. For $n=2$, we show that $p(x,y,x_1,x_2)=(x+y)(yx+(x+x_1+y)(x+x_2+y))$ is Lorentzian, from which we conclude that $\homo {L(\L_2, x, x_1,x_2)} = p(x,y,x_1,x_2) \cdot (x+y) $ is Lorentzian. The support satisfies $\supp((x+y)(yx+(x+x_1+y)(x+x_2+y))) = \supp((x+y)(x+x_1+y)(x+x_2+y))$, which is the support of a Lorentzian polynomial, thus M-convex. We proceed with checking the eigenvalue condition. If we take a partial derivative with respect to $x_1$ or $x_2$ we end up with a Lorentzian polynomial, since $\partial_{x_1} p(x,y,x_1,x_2) = (x+y)(x+x_2+y)$ is a product of Lorentzian polynomials (and similarly for $\partial_{x_2})$. Thus, it is enough to show that the Hessians of $\partial_x  p(x,y,x_1,x_2)$ and $ \partial_y p(x,y,x_1,x_2)$ have at most one positive eigenvalue. It suffices to check the Hessian of $\partial_x  p(x,y,x_1,x_2)$, because $p$ is symmetric in $x$ and $y$. It has the form 
    \begin{align*}
      H =\begin{pmatrix}
           0& 1& 2 & 2 \\
           1& 0& 2& 2   \\
           2& 2& 6& 8  \\
           2& 2& 8& 8
         \end{pmatrix}, 
    \end{align*}
    and the eigenvalues are $\lambda_1 = \frac{3}{2}(5+\sqrt{33}), \lambda_2 = - \frac{3}{2}(\sqrt{33}-5), \lambda_3 = -1, \lambda_4 =0$.  
    
    Let $n$ be greater than 2. As $(x+y)^n$ is Lorentzian, the result follows by proving  
    \begin{align*}
      p(x,y,(x_i)) = xy^{n-1} + \prod_{i=1}^{n}(x+x_i+y)
    \end{align*}
    to be Lorentzian after being multiplied by $(xy)^k$, for some $k \in \bN$ (note that this does not hold when $n=2$). We start by proving M-convexity. This follows directly from the fact that 
    \begin{align*}
      \supp(p) = \supp\left(xy^{n-1} + \prod_{i=1}^{n}(x+x_i+y)\right) = \supp\left(\prod_{i=1}^n (x+x_i+y)\right),
    \end{align*}
    as $\prod_{i=1}^n (x+x_i+y)$ is a product of Lorentzian polynomials and hence M-convex (recall also the discussion in Definition~\ref{def:pre-lorentzian-tree} which established that the $(xy)^k$ factor may be ignored when checking M-convexity).

    It remains to be shown that for every $(n+2k-2)$th partial derivative of
    \begin{align}
      \begin{split} \label{eq:q-poly}
        q(x,y,(x_i)) &\defeq (xy)^k\cdot \left(xy^{n-1} + \prod_{i=1}^{n}(x+x_i+y)\right)\\
                     &= x^{k+1}y^{n+k-1} + x^ky^k\prod_{i=1}^{n}(x+x_i+y),                     
      \end{split}
    \end{align}
    there exists $k \in \bN$ such that the associated Hessian has exactly 1 positive eigenvalue. In fact, we will show that this holds for all $k \in \bN$ sufficiently large. We note that, for any $j \in [n]$, we have
    \begin{align*}
      \partial_{x_{j}}q(x,y,(x_i)) = x^k y^k\prod_{i \neq j}(x+x_i+y), 
    \end{align*}
    which is a product of Lorentzian polynomials and therefore Lorentzian. We can therefore restrict ourselves to checking that partial derivatives of the form
    \begin{align*}
      \partial_x^l \partial_y^{n+2k-l-2}q(x,y,(x_i))
    \end{align*}
    have Hessians with exactly one positive eigenvalue (for large enough $k$). We now proceed by splitting this problem in four cases based upon the value of $l$.

    \begin{mycase}{1: $l < k-1$ or $l > k+1$}
      This is the simplest case as such sequences of partial derivatives all kill the $y^{n+k-1}x^{k+1}$ term appearing in~\eqref{eq:q-poly}, leaving just $x^ky^k\prod_{i=1}^{n}(x+x_i+y)$ which is a product of Lorentzian polynomials and hence Lorentzian.
    \end{mycase}
    
    In order to tackle the remaining cases, we will need to make use of Lemma~\ref{lem:lemma1}, so let us now justify its relevance. As $q(x,y,(x_i))$ is symmetric in the $x_i$ variables and $\deg_{x_i}(q)=1$, the Hessian associated to any $(n+2k-2)$th partial derivative of $q$ satisfies the conditions of Lemma~\ref{lem:lemma1}.

    \begin{mycase}{2 : $l = k-1$} In this case, the only monomials in $q$ which survive and contribute to the Hessian are the ones that have degree in $x$ at least $k-1$ and degree in $y$ at least $n+k-1$. In particular, they are given by
      \begin{gather*}
        x^{k-1}y^{n+k-1}x_ix_j,\, x^{k}y^{n+k-1}x_i,\, x^{k-1}y^{n+k}x_i, \\ x^{k+1}y^{n+k-1},\, x^{k-1}y^{n+k+1} \ \text{and} \  x^{k}y^{n+k}.
      \end{gather*}
      Expanding $q$, we see that their respective coefficients are 
      \begin{align*}
        0,\, 1,\, 0,\, n+1,\, 0 \ \text{and} \  1.
      \end{align*}
      After differentiating and removing the common $(k-1)!(n+k-1)!$ factor, which comes from the exponents in $x$ and $y$, we are left with a Hessian $H$, whose reduced form (see Lemma~\ref{lem:lemma1}) is given by
      \begin{align*}
        \red{H} =
        \begin{pmatrix}
          0& k & 0 \\
          nk& k(n+1)(k+1)& k(n+k)  \\
          0& k(n+k) & 0
        \end{pmatrix}.
      \end{align*}
      Developing the first row gives $\det \red{H} = 0$. Furthermore, the trace is positive and    
      \begin{align*}
        \red{H}_{11} + \red{H}_{22} + \red{H}_{33} = -nk^2 - k^2(n+k)^2 < 0.
      \end{align*}
      We can therefore apply Lemmas~\ref{lem:lemma1} and~\ref{lem:lemma2} to conclude that $H$ has exactly $1$ positive eigenvalue. 
    \end{mycase}

    We note that up to this point in the proof, the $(xy)^k$ factor (which, we recall, characterises the difference between $p$ being pre-Lorentzian and $\homo{p}$ simply being Lorentzian) has not played a role. However, this will not be the case for the remaining cases, which will require us to study the asymptotics in terms of $k$.

    \begin{mycase}{3 : $l = k$} In this case, the only monomials in $q$ which survive and contribute to the Hessian are the ones that have degree in $x$ at least $k$ and degree in $y$ at least $n+k-2$. In particular, they are given by
      \begin{gather*}
        x^ky^{n+k-2}x_ix_j,\, x^{k+1}y^{n+k-2}x_i,\, x^ky^{n+k-1}x_i, \\ x^{k+2}y^{n+k -2},\, x^ky^{n+k} \ \text{and} \  x^{k+1}y^{n+k -1},
      \end{gather*}
      with coefficients in $q$ given by
      \begin{align*}
        1,\, n-1,\, 1,\, \binom{n }{2},\, 1 \ \text{and} \ n+1,
      \end{align*}
      respectively. As before, after taking derivatives we can factor out $k!(n+k-2)!$ and we get the reduced Hessian
      \begin{align*}
        \red{H}=
        \begin{pmatrix}
          n-1& (k+1)(n-1)& n+k-1 \\
          n(k+1)(n-1)& (k+2)(k+1)\binom{n }{2}&(k+1)(n+k-1)(n+1)    \\
          n(n+k-1)&(k+1)(n+k-1)(n+1)  & (n+k)(n+k-1)
        \end{pmatrix}.
      \end{align*}
      Note that, for arbitrary $n \in \bN$ and $k \leq 2$, this matrix has two positive eigenvalues. However, we claim that for $k$ large enough this cannot happen. In fact, as before, this will follow from Lemma~\ref{lem:lemma2}. As the trace is still clearly positive we just have to show that $\red{H}_{11} + \red{H}_{22} + \red{H}_{33} \leq 0$ and $\det(\red{H}) \geq 0$. To see this, first note that $\red{H}_{22}$ and $\red{H}_{33}$ are both of order $k^2$, whereas 
      \begin{align*}
        \red{H}_{11} = k^4\left(\binom{n }{2}-(n+1)^2\right) + O(k^3),
      \end{align*}
      which is clearly negative for large enough $k$. For the determinant, we consider the following simplified matrix $\red{H}'$ where we remove any factors containing $k$:
      \begin{align*}
        \red{H}' =\begin{pmatrix}
                    n-1& n-1& 1 \\
                    n(n-1)& \binom{n }{2}&n+1   \\
                    n&  n+1  & 1
                  \end{pmatrix}.
      \end{align*}
      In particular, expanding $\det(\red{H})$ gives
      \begin{align*}
        \det(\red{H}) = k^4\cdot \det(\red{H}') + O(k^3),
      \end{align*}
      implying that it is enough to see that $\det(\red{H}')>0$. A simple calculation gives
      \begin{align*}
        \det(\red{H}') = \frac{n^2-3n}{2}+1,
      \end{align*}
      which is positive for $n>2$. 
    \end{mycase}

    \begin{mycase}{4 : $l = k+1$} In this case, the only monomials in $q$ which survive and contribute to the Hessian are the ones that have degree in $x$ at least $k+1$ and degree in $y$ at least $n+k-3$. In particular, they are given by
      \begin{gather*}
        x^{k+1}y^{n+k-3}x_ix_j,\, x^{k+2}y^{n+k-3}x_i,\, x^{k+1}y^{n+k-2}x_i,\\ x^{k+3}y^{n+k-3},\, x^{k+1}y^{n+k-1} \ \text{and} \ x^{k+2}y^{n+k-2},
      \end{gather*}
      with coefficients in $q$ given by
      \begin{align*}
        n-2,\, \binom{n-1 }{2},\, n-1,\, \binom{n }{3},\, n+1 \ \text{and} \ \binom{n }{2}
      \end{align*}
      respectively. As before, after taking derivatives we can factor out $(k+1)!(n+k-3)!$ and we get the reduced Hessian
      \begin{align*}
        \red{H}=
        \begin{pmatrix}
          (n-1)(n-2)& (k+2)\binom{n-1}{2}& (n+k-2)(n-1) \\
          n(k+2)\binom{n-1}{2}& (k+3)(k+2)\binom{n}{3}&(k+2)(n+k-2)\binom{n}{2}    \\
          n(n+k-2)(n-1) &(k+2)(n+k-2)\binom{n}{2}  & (n+1)(n+k-1)(n+k-2)
        \end{pmatrix}.
      \end{align*}
      Again, the trace is clearly positive, and we only have to check that $\red{H}_{11} + \red{H}_{22} + \red{H}_{33} \leq 0$ and $\det(\red{H}) \geq 0$. The argument is extremely similar to the one given in Case 3: $\red{H}_{22}$ and $\red{H}_{33}$ are both of order $k^2$, whereas 
      \begin{align*}
        \red{H}_{11} = k^4\left(\binom{n }{3}(n+1) - \binom{n }{2}^2\right) + O(k^3), 
      \end{align*}
      which is clearly negative for large enough $k$. Finally, for the determinant, we again simplify the matrix by removing the factor with any $k$ dependence:
      \begin{align*}
        \red{H}' =
        \begin{pmatrix}
          (n-1)(n-2)& \binom{n-1 }{2}& n-1 \\
          n\binom{n-1 }{2}& \binom{n }{3}&\binom{n }{2}   \\
          n(n-1)&  \binom{n }{2}  & (n+1)
        \end{pmatrix},
      \end{align*}
      and compute
      \begin{align*}
        \det(\tilde{H}') = \frac{n^4 - 4n^3 + 5n^2 - 2n}{6},
      \end{align*}
      which is positive for $n>2$. \endproof
    \end{mycase}
  \end{THM}

  \begin{REM}
    \label{rem:Ln-hypotheses}
    We now make a few remarks about the hypotheses of Theorem~\ref{thm:leafy-stars-lorentz}. One might ask whether $C(\L_{n},x,(x_i))$ is real stable, which is a glueable property. This is false, since $I(\L_3;x)$ is not real rooted and this approach already fails for $n=3$. We also observe two interesting features of $\L_n$. If we remove all leaves from $\L_n$, the support of the coloured star we get is not M-convex and a naive attempt of glueing stars fails as well. If we instead only remove the leaves attached to the center of $\L_n$ we indeed end up with a support that is M-convex, but the eigenvalue condition fails. So, in some sense the additional leaves attached to the center \say{push} the eigenvalues of the Hessian in just the right way to guarantee that the graph is pre-Lorentzian. 
  \end{REM}

  We can now apply the glueing procedure established in Proposition~\ref{prop:pre-lorentzian-glue} for the following result. 

  \pasttheorem{thm:main-theorem}
  \proof 
  Since being pre-Lorentzian is stable under disjoint union (indeed, this follows from Proposition~\ref{prop:lorentz-polyn-stable-indent}, Proposition~\ref{prop:lorentz-polyn-stable-prod} and the fact that $\L(\G_1 \sqcup \G_2) = \L(\G_1)\cdot \L(\G_2)$), in light of Lemma~\ref{lem:pre-lorentzian-log-concave}, the result follows from the claim that every connected graph in $R_{W_4}(\graphs)$ is pre-Lorentzian.

  Let $G$ be connected in $R_{W_4}$. Proposition~\ref{prop:replace_to_glue} then tells us that $G$ may be constructed by glueing copies of $\L_n$ across free vertices. By Theorem~\ref{thm:leafy-stars-lorentz} these copies are all pre-Lorentzian, and furthermore, by Proposition~\ref{prop:pre-lorentzian-glue}, so is $G$ (thought of as a partitioned graph with no free vertices). 
  \endproof 

  \section{Conclusion}

  While our main result, Theorem~\ref{thm:main-theorem}, does, in a sense, fully utilise the content of Theorem~\ref{thm:leafy-stars-lorentz} (as Proposition~\ref{prop:replace_to_glue} in fact describes a bijection between $R_{W_4}$ and the set of glueings of $\L_n$ with no remaining free variables), there are broader ways in which the results of this paper can be extended. We shall now briefly describe two such directions.

  Empirical evidence suggests that results analogous to Theorem~\ref{thm:leafy-stars-lorentz} hold for other families of labelled graphs. Such results could be combined with Theorem~\ref{thm:leafy-stars-lorentz} to prove that arbitrary glueings of graphs from any such families are pre-Lorentzian and hence have log-concave independence sequences, broadening the scope of Theorem~\ref{thm:main-theorem}.  

  Alternatively, one could try to adapt our approach to directly attack Alavi, Malde, Schwenk and Erd\H{o}s' conjecture that all forests have a unimodal independence sequence. In particular, one could further weaken the pre-Lorentzian condition, aiming to find a criteria which directly implies unimodality of the independence sequence (and not log-concavity, circumnavigating Theorem~\ref{thm:trees-not-log-concave}) all while preserving its good behaviour with respect to the glueing of labelled graphs.   

  \bibliographystyle{alpha}
  \bibliography{biblio}

\end{document}